\documentclass[runningheads]{llncs}
\usepackage[T1]{fontenc}
\usepackage{graphicx}
\usepackage{amssymb} 
\usepackage{mathtools}
\usepackage{soul} 
\usepackage{xcolor} 

\usepackage{amsmath}
\usepackage{multirow}
\usepackage{tabularx}
\usepackage{booktabs}

\begin{document}
\title{Learning to Choose Branching Rules for Nonconvex MINLPs}
%
%
\author{Timo Berthold\inst{1}\orcidID{0000-0002-6320-8154} \and
Fritz Geis\inst{2}}
\authorrunning{T.~Berthold and F.~Geis}
%
\institute{Fair Isaac Germany GmbH, Takustr. 7, 14195 Berlin, Germany
\email{timoberthold@fico.com} \and
Zuse Institute Berlin, Takustr. 7, 14195 Berlin, Germany\\
\email{fritzgeis@gmail.com}
}
\maketitle              

\begin{abstract}
Outer-approximation-based branch-and-bound is a common algorithmic framework for solving MINLPs (mixed-integer nonlinear programs) to global optimality, with branching variable selection critically influencing overall performance. In modern global MINLP solvers, it is unclear whether branching on fractional integer variables should be prioritized over spatial branching on variables, potentially continuous, that show constraint violations, with different solvers following different defaults. We address this question using a data-driven approach. Based on a test set of hundreds of heterogeneous public and industrial MINLP instances, we train linear and random forest regression models to predict the relative speedup of the FICO\textsuperscript{®} Xpress Global solver when using a branching rule that always prioritizes variables with violated integralities versus a mixed rule, allowing for early spatial branches.

We introduce a practical evaluation methodology that measures the effect of the learned model directly in terms of the shifted geometric mean runtime. Using only four features derived from strong branching and the nonlinear structure, our linear regression model achieves an 8--9\% reduction in geometric-mean solving time for the Xpress solver, with over 10\% improvement on hard instances.
We also analyze a random regression forest model. Experiments across solver versions show that a model trained on Xpress 9.6 still yields significant improvements on Xpress 9.8 without retraining.

Our results demonstrate how regression models can successfully guide the branching-rule selection and improve the performance of a state-of-the-art commercial MINLP solver.

\keywords{Nonlinear Optimization  \and Machine Learning \and Branching.}
\end{abstract}

\section{Introduction}\label{sec:Intro}

We consider \emph{MINLPs} \emph{(mixed-integer nonlinear programs)} of the form 
\begin{align}
    min\{c^T x\mid g_k(x)\leq 0,\forall k\in\mathcal{K}, l\leq x\leq u, x_j\in\mathbb{Z}, \forall j\in\mathcal{J}\},
    \label{minlp_def}
\end{align} where all constraint functions $g_k:\mathbb{R}^{n}\rightarrow\mathbb{R}$ are factorable and all variable bounds $l,u\in\bar{\mathbb{R}} \coloneq\mathbb{R}\cup\{\pm \infty\}$. The set $\mathcal{K}=\{1,\dots,m\},m\in\mathbb{N}$, indexes the constraints and $\mathcal{J}\subseteq\{1,\dots, n\}$ the integer variables. A nonlinear objective can be easily modeled by introducing an auxiliary variable and an objective-transfer constraint, see, e.g.,~\cite{Vigerske2013}. If all $g_k$ are linear, and $\mathcal{J}=\emptyset$ we call (\ref{minlp_def}) a \emph{linear program (LP)}. This work focuses on \emph{nonconvex MINLPs}, i.e., problems of form~(\ref{minlp_def}), where at least one $g_k$ is nonconvex.

Note that factorable functions can be represented via a directed acyclic \emph{expression graph}, with nodes representing operators or variables, and arcs representing the data flow of the computation. In this paper, we refer to this representation as the \emph{DAG}, for a good overview on the use of the DAG in MINLP solving, we recommend \cite{Vigerske2013}.

For solving problems of the form (\ref{minlp_def}), we use the  FICO\textsuperscript{\textregistered} Xpress Global \cite{xpress} MINLP solver, which we will refer to as \emph{Xpress}. 
Xpress is based on the \emph{branch-and-bound} method (\emph{B}\&\emph{B}), which recursively partitions the problem by splitting the domain of selected variables, which is called \emph{branching}. Selecting good branching variables is crucial for the performance of B\&B-based MINLP solvers, see, e.g.,~\cite{BelottiKirchesLeyfferLinderothLuedtkeMahajan2013}. For more details on the implementation in Xpress, see \cite{BelottiBertholdGallyGottwaldPolik2025}. 

This paper studies a fundamental question:  should we always branch on fractional integer variables first, or consider spatial branches even when there are fractional integers? Unlike most prior work on using ML for branching \cite{khalil2016,scavuzzo2024,nair2020,gupta2020,alvarez2017}, we do not learn individual branching decisions or attempt to mimic existing strategies such as strong branching; instead, we perform algorithm selection between two established branching rules. Different from most prior work, we consider a heterogeneous set of instances. The resulting model can be integrated directly into solver code and does not require any pre-training on the user side. This is akin to prior ML-based algorithm-selection work to choose between scaling procedures~\cite{BertholdHendel2021}. local cut selection rules~\cite{BertholdFrancobaldiHendel2025}, linearization techniques~\cite{bonami2022classifier}, or spatial branching strategies in the context of RLT for polynomial optimization~\cite{ghaddar2023learning}, respectively, and deliberately different from solver-free learning approaches for MINLP as in~\cite{tang2024learning}.

\section{A Quick Recap of Branching for MINLPs}
The B\&B algorithm recursively partitions the problem into smaller subproblems (\emph{branching}) and solves LP relaxations to obtain bounds (\emph{bounding}) until an optimal solution or infeasibility proof is found.

An \emph{LP relaxation} of an MINLP
is obtained by dropping integrality constraints and replacing nonlinear constraints with linear underestimators where possible. 
This relaxation is successively strengthened by \emph{outer-approximation cuts}  \cite{DuranGrossmann1986}. A well-designed cutting plane separation procedure often helps to reduce the branch-and-bound tree size while accelerating the overall
solving process  \cite{TurnerBertholdBesanconKoch2023}. Unlike in MIP solving, cutting planes are often additionally separated immediately during branching-node creation in MINLPs.

In this paper, we focus on \emph{variable branching}, in which the domain of a single variable is split into two or more intervals.

Two key types of variable branching are:
\begin{enumerate}
    \item \emph{Integer branching}, which is applied when an integer variable has a fractional value $\check{x}_j\in\mathbb{R}\setminus\mathbb{Z}$ in the solution $\check{x}$ of the current LP relaxation. Two subproblems are created that enforce $x_j\leq \lfloor \check{x}_j\rfloor$ and $x_j\geq \lceil \check{x}_j\rceil$, respectively.
    \item \emph{Spatial branching} \cite{HorstTuy1996}
is applied when the violation of a nonconvex constraint cannot be resolved by an outer-approximation cut, but requires partitioning variable domains. Spatial branching candidates are often continuous variables, but can also include integer variables whose LP value happens to be integral. Two created subproblems enforce
$x_j\leq \lfloor \check{x}_j\rfloor$ and $x_j\geq \lceil \check{x}_j\rceil$, respectively, for a branching point $\check{x}_j\in\mathbb{R}$. Though the LP solution is not explicitly excluded, subsequent outer-approximation cuts typically remove it.
\end{enumerate}

\section{Machine Learning Methodological Approach}\label{sec:ML}

\subsubsection{Learning Task/Feature Space} Our learning task consists of choosing, after root node processing and right before the first branch, one of two rules of how to combine integer branching and spatial branching for the remainder of the branch-and-bound search: Either, always branch on integer candidates and conduct spatial branches only when there is no integer branching candidate, which we will refer to as "PreferInt" (this is the default, e.g., in the SCIP MINLP solver). Or mix both candidate sets and always allow the choice of either type of candidate (which is the default, e.g., in the Xpress solver), which we refer to as "Mixed".\footnote{We ruled out always preferring spatial branches in a preliminary experiment, since this option was a factor eight slower on average and rarely won against the others.} 

Although this is inherently a binary decision, we frame it as a regression problem. This choice is motivated by two considerations. Firstly, our ultimate goal is to improve the average runtime of the solver, which is a metric that is numerical and not categorical. Secondly, our focus is on getting the prediction right for those instances on which the performance of selecting "Mixed" and "PreferInt" significantly differs, see also  \cite{BertholdFrancobaldiHendel2025}. Regression allows us to model the magnitude of this difference directly and thereby focus the learning on the cases where it matters most.

To this end, we train regression models  $y_i:\mathbb{R}^{d}\rightarrow\mathbb{R}$ that map a $d$-dimensional feature vector $f=(f_1, \dots, f_d)$  onto the speedup or slowdown factor (the \emph{label}) in runtime by using "PreferInt" instead of "Mixed".
We initially used 17 features, see Table \ref{tab:FeatureDef}.

This includes features related to strong branching at the root node, such as the average change in the dual bound resulting from integer and spatial strong branching, \texttt{AvgRelBndChngSBLPInt} and \texttt{AvgRelBndChngSBLPSpat}, respectively,
the number of variables fixed from strong branching on spatial branching candidates \texttt{\#SpatBranchEntFixed}
\footnote{There were only a few instances where integer strong branching fixed variables; hence, a corresponding feature would have been almost flat zero.},
and the amount of deterministic \emph{work} invested in either strong branching, \texttt{AvgWorkSBLPInt} and  \texttt{AvgWorkSBLPSpat}.  Work is a deterministic measure of computational effort implemented in 
Xpress. These features give us an indication of how effective (and expensive) strong branching on integer or spatial variables is. Relatedly, \texttt{\#IntViols} and \texttt{\#NonlinViols} refer to the number of integer and spatial branching candidates. 

As problem structure features, we include the percentage of variables that are integer, \texttt{\%IntVars}, the percentage of constraints that are equations, \texttt{\%EqCons}, the ratio of quadratic elements in the problem to variables, \texttt{\%QuadrElements},  and the percentage of constraints that contain nonlinearities, \texttt{\%NonlinCons}.
Further, to measure the nonlinearity of the problem, we use information about the DAG, in particular, the percentage of variables that are part of any nonlinearity, \texttt{\%VarsDAG}, and the ratio between nodes in the DAG and nonzeros in the linear part of the problem, \texttt{NodesInDAG}. We further include the percentage of integer and unbounded variables among all variables in the DAG, \texttt{\%VarsDAGInt} and \texttt{\%VarsDAGUnbnd}, 
as for integer DAG variables, we "hit two birds with one stone" and branching on unbounded variables can be crucial to get efficient dual bounds. Finally, we consider \texttt{\%QuadrNodesDAG} to measure whether the nonlinearities in the problem are mostly quadratic.


\begin{table}[!ht]
\renewcommand{\arraystretch}{1.1}
\centering
\resizebox{0.9\textwidth}{!}{%
\begin{tabularx}{\textwidth}{|l X|}
    \hline
     \textbf{Feature} & \textbf{Feature Scaling} \\
    \hline
 \multicolumn{2}{|c|}{\emph{Problem Structure}} \\
    \hline
     \texttt{\%QuadrElements} & number quadratic elements over $n$ \\
     \texttt{\%IntVars} & \#Integer variables after presolve over $\tilde{n}$ \\
     \texttt{\%EqCons} & \#equality constraints over $m$ \\
     \textbf{\texttt{\%NonlinCons}} & \#nonlinear constraints over $m$ \\
    \hline
     \multicolumn{2}{|c|}{\emph{Effect of Branching}} \\
    \hline
     \texttt{\#IntViols} & \\
     \texttt{\#NonlinViols} & \\
     \texttt{\#SpatBranchEntFixed} & \\
     \texttt{AvgWorkSBLPInt} & \\
     \texttt{AvgWorkSBLPSpat} & $\log_{10}(\text{Value+1})$ \\
     \texttt{AvgRelBndChngSBLPInt} & \\
     \texttt{AvgRelBndChngSBLPSpat} & \\
     \texttt{AvgCoeffSpreadConvCuts} & \\
    \hline
     \multicolumn{2}{|c|}{\emph{DAG}} \\
    \hline
     \texttt{NodesInDAG} & NodesInDAG over NodesInDAG+$\tilde{M}$ \\
     \texttt{\%VarsDAG} & \#vars in DAG over $\tilde{n}$ \\
     \texttt{\%VarsDAGUnbnd} & \#unbounded vars over \#vars in DAG \\
     \texttt{\%VarsDAGInt} & \#integer vars over \#vars in DAG \\
     \texttt{\%QuadrNodesDAG} & \#quadratic operator nodes in DAG over all\\ 
     & nonlinear operator nodes in DAG \\
    \hline
\end{tabularx}
}
\caption{$m$ and $n$ are the number of constraints and variables before presolving, respectively; $\tilde{n}$ and $\tilde{M}$ the number of variables and linear nonzeros after presolving.}
\label{tab:FeatureDef}
\end{table}

\subsubsection{Data}
The data on which the models are trained comes from running Xpress 9.6\footnote{More precisely: An internal version of the Xpress 9.6 that exposes those features that are otherwise not available as public attributes.} twice on a heterogeneous benchmark of 683 public and industrial MINLP instances, each with two permutations to mitigate the effect of performance variability~\cite{lodi2013performance,gamrath2020exploratory}, yielding 2049 data points. For each instance, we record the runtimes produced by both branching rules and the complete feature set. Instances solved at the root or otherwise unsuitable for comparison are filtered out, resulting in a final dataset of 797 data points, see  \cite{Geis2025} for details. Solving at the root node was by far the most common reason for filtering.


\subsubsection{Training} For training the models, we split the data randomly into $80\%$ training and $20\%$ test set. 
The models we train on the training set are a linear regressor  \cite{ChambersHastie1992} and a random forest regressor, \texttt{RF},  \cite{LiawWiener2002}. We use the python library scikit-learn  \cite{scikit-learn}, which provides us with the linear regressor by the function \emph{LinearRegression} and the random forest regressor by the function \emph{RandomForestRegressor}.

\subsubsection{Testing} Instead of training one linear regressor and one random forest regressor, we opted for training and testing one hundred models each with different random seeds and average their performance scores to evaluate how promising this ML-based approach is. 

\label{sssec:Measure}
 To measure the performance of the regression models, we use the \emph{accuracy} and the shifted geometric mean of the runtime (\emph{sgm\_runtime}). The accuracy is defined as the percentage of times the model predicted the faster rule. The sgm\_runtime is the shifted geometric mean of runtimes when solving each test instance using the predicted branching rule over the shifted geometric mean time using always the default rule.
  
 Hence, accuracy is always between 0 and 100\%, with larger values being better. Sgm\_runtime can be smaller or larger than 1, with values larger than 1 indicating a deterioration and values smaller than 1 indicating an improvement: the smaller the number, the better. This is the primary performance indicator for solver development in practice. 
 
 To compute the shifted geometric mean  \cite{Achterberg2007} with a shift of $10$,  
 measurements $X = (X_1,\dots,X_n)$ are aggregated via
$sgm(X)=-10+\prod\limits_{i=1}^{n}(X_i+10)^{\frac{1}{n}}$. The use of the shifted geometric mean is a commonly used method to aggregate performance measures, in particular running time, inn mathematical optimization~\cite{BertholdHendel2021}.

For the linear model, \emph{feature importance} is given by the absolute value of the learned coefficients, whereas for the random forest it is measured as the normalized total reduction in mean squared error induced by splits on that feature (mean decrease impurity, MDI), which are the default importance metrics in scikit-learn.
For each random seed and each type of model, we computed the feature importance for all features and sorted them from most important to least important. Then we assigned a score of zero to the most important one, a score of one to the second most important one, and so on. Finally, we added, for each model type, the scores across all one hundred runs together. The four most important features per model type (as by this score sum) are listed in Table \ref{tab:FeatImpo}. 
Although the top-ranked features differ between the linear regression and random forest models, there is overlap in terms of the underlying information captured. 
In particular, \texttt{AvgRelBndChngSBLPSpat} ranks first for the linear model and second for the random forest, and \texttt{\%NonlinCons}, ranked third for the linear model, is a very close fifth for the random forest.
Overall, five of the eight highest-ranked features coincide across the two models. 
Differences are expected given the different nature of the models: linear regression emphasizes globally predictive, approximately linear effects with low collinearity, whereas random forests prioritize features that enable strong local splits and nonlinear interactions, for instance, capturing cases where a branching rule is beneficial for either extreme but not for intermediate feature values.

\subsubsection{Further Approaches} 
In earlier versions, we tested the algorithm with different features, unscaled or differently scaled features,
 on an earlier version of Xpress and for the SCIP solver (where it also improved performance, but not as much as in the Xpress case).
Details can be found in the thesis \cite{Geis2025}.
This thesis also contains a detailed description of how we selected the scaler and imputer for the data set and a discussion of restricting the decision tree depth to five in the random forest models. 
\begin{table}
    \centering
    \begin{tabular}{|c|l|l|}
        \hline
         \textbf{Ranking} & \textbf{Linear}& \textbf{Forest} \\
        \hline
          1. & \texttt{AvgRelBndChngSBLPSpat} & \texttt{AvgCoeffSpreadConvCuts} \\
          2. & \texttt{\%IntVars}        & \texttt{AvgRelBndChngSBLPSpat} \\
          3. & \texttt{\%NonlinCons}            & \texttt{\#NonlinViols} \\
          4. & \texttt{\%VarsDAGInt}          & \texttt{\%EqCons} \\
        \hline
    \end{tabular}
    \caption{Four most important features for either model type.}
    \label{tab:FeatImpo}
\end{table}
\section{Computational Experiments}
Our computational experiments consist of three parts: 
Firstly, we evaluate the regression models trained on the full 17-feature set, and then examine how their performance evolves as we iteratively remove the least important features. This reduction process provides insights into which features drive prediction quality and whether a more compact feature set can yield comparable performance with presumably better robustness.
Secondly, we analyze the final reduced models in more detail
and provide an analysis of their performance with respect to accuracy and runtime.
Finally, we compare how those models continue to perform as the underlying branch-and-bound method is improved (in this case, through a solver version update).

\begin{figure}[ht!]
    \centering
    \begin{minipage}[t]{0.48\textwidth} 
        \centering
    \includegraphics[width=1\linewidth]{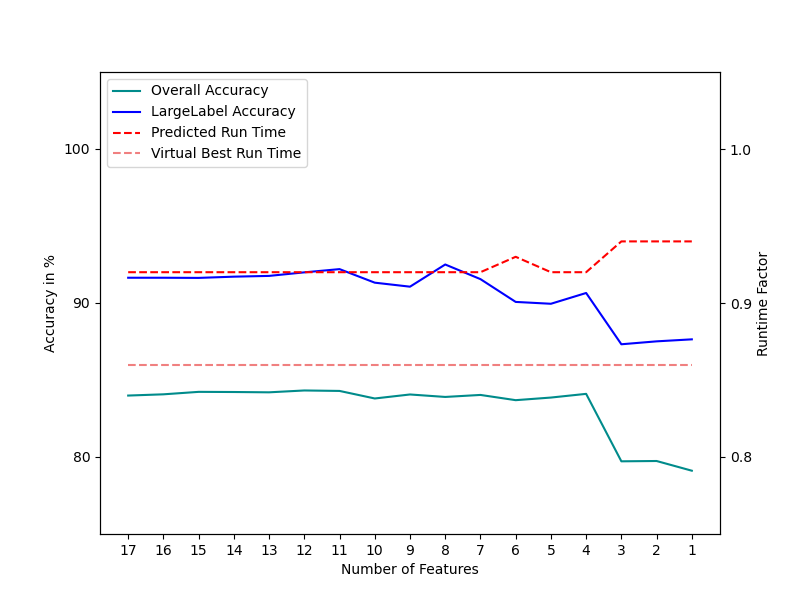}
    \caption{Impact of feature-reduction on accuracy and runtime for linear regression models}
    \label{fig:96_Lin_Reduction}
    \end{minipage}
    \hfill
    \begin{minipage}[t]{0.48\textwidth} 
    \centering
        \includegraphics[width=\linewidth]{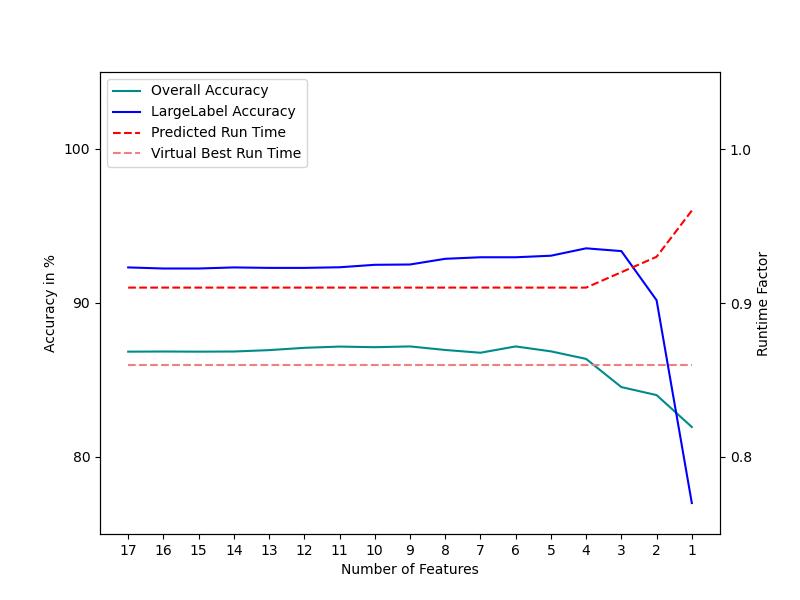}
        \caption{Impact of feature-reduction on accuracy and runtime for random forest regression models}
        \label{fig:96_For_Reduction}
    \end{minipage}
\end{figure}

For feature-reduction in the first experiment, we use
 the ranking by 100 different seeds described in Section \ref{sssec:Measure}. We remove one feature at a time (the least important) and retrain after each step.
The impact on both accuracy and sgm\_runtime for the linear models and the random forest models is shown in Figures \ref{fig:96_Lin_Reduction} and \ref{fig:96_For_Reduction}, respectively. The $x$-axis shows the number of features used, the $y$-axis shows the two performance measures. The dotted red line represents the sgm\_runtime, and the dotted salmon-colored constant represents the best possible sgm\_runtime score, achievable only if a model made a correct prediction in all cases. The green solid line represents the accuracy score on all test instances, while the blue shows the accuracy score on instances with a label, i.e., a speedup or slowdown factor, larger than four, here called \emph{LargeLabel Accuracy}.
These are the most important instances to predict accurately when the goal is to reduce the mean runtime. 

In Figure \ref{fig:96_Lin_Reduction} we can see that the performance of the initial linear model with all features is at about 85\% overall accuracy, more than 90\% LargeLabel accuracy, and a sgm\_runtime factor of a little over 0.91. The accuracy and the sgm\_runtime factor of the random forest is a bit better, at 87\% and a little over 0.9, respectively. LargeLabel accuracy is about the same.
The feature-reduction experiments reveal a remarkably stable behavior: both model types maintain essentially constant accuracy and runtime performance as long as at least four features are retained. The LargeLabel accuracy drops slightly when going down from 17 to 4 features.
Once the feature count drops below this threshold, all measures deteriorate noticeably. 
Based on the aggregated importance scores, we identified the four most influential features for each model type, see Table~\ref{tab:FeatImpo}, and use them for the subsequent experiments.

Training and testing both model types using only their respective four most important features yields the performance scores depicted in Table \ref{tab:MainResult}. The scores are divided by model type and, for checking possible overfitting, by training and test set. Again, the accuracy score is calculated \emph{Overall}  and on LargeLabel instances. The two bottom rows contain the sgm\_runtime scores. As before, all scores are averages across one hundred runs for each model type.

 The linear model performs almost identical on training and test set; there is no indication of overfitting. On the test set, it predicts the better branching rule with an accuracy of about 84\% overall and more than 90\% on the subset of instances with the largest performance differences. Most importantly, the sgm\_runtime score is $0.919$ on the test set, indicating a significant 8\% speedup. The standard deviation was 2.8\%. On models that take more than 100 seconds to solve by either of the two variants, we observed a speedup of over 10\% with a standard deviation of 6.8\%. 
  These results show that a linear model, based on a handful of carefully designed features, provides meaningful and robust performance gains, with an easy-to-implement and interpretable model, which is highly desirable in practice.

For the random forest, however, there is a significant difference between the training and the test set. The almost perfect accuracy on the training set drops to $93.6\%$ on the test set, a clear sign of overfitting, and the sgm\_runtime clearly deteriorates as well. With a value of 3.2\%, the standard deviation was also higher compared to the linear model. From a practical perspective, the linear model seems preferable, even though its performance is slightly worse even on the test set.

The linear model predicted the "Mixed" rule to be faster for around 60\% of the instances and "PreferInt" to be faster for 40\% of the instances. For the random forest, the split is close to 50:50. 
For the linear model, in 51\% of the cases, switching to the "PreferInt" rule improves performance by at least 10\%, in 11\% of the cases performance got at least 10\% worse, in the remaining 38\%, it stayed roughly the same. For the random forest, this split was 43\% wins, 11\% losses, 36\% within $\pm$10\% runtime.


With our final experiment, we approach the question of how sensitive the results are to the changes in the underlying MINLP solving algorithm, more specifically, the robustness across solver versions. Therefore, we apply the models trained on Xpress 9.6 to data generated by Xpress 9.8. The latter version introduces substantial algorithmic improvements that make it roughly 50\% faster on average for difficult MINLP instances, including changes to presolving, cut generation, heuristics, and branching logic. Despite these differences, both learned models continue to yield computational benefits: the linear model achieves an average speedup of about
$3.3\%$ and the random forest about $4.5\%$.


This leaves the question of whether the smaller performance gain is inherent to the other changes in the solver or to the trained models not being as accurate for predicting performance of the newer solver version. In an additional experiment, we trained new models from scratch on the data of Xpress 9.8, and to our positive surprise, both the accuracy and the time factors were very similar to the ones shown in Table \ref{tab:9698}. The linear model trained on Xpress 9.8 data achieved the same sgm\_runtime score of $0.967$, whereas the random forest slightly improved, presumably again a result of overfitting. Furthermore, the sets of the four most important features were almost identical, with only one feature differing in the linear model and none in the random forest. 
This robustness is particularly valuable for practical deployment, as it suggests that models need not be retrained for every release cycle.

\begin{table}[h!]

\begin{minipage}{0.5\textwidth}
\centering
\begin{tabular}{|c|c|c|c|c|}
        \hline
        & \multicolumn{2}{c|}{\textbf{Linear}} & \multicolumn{2}{c|}{\textbf{Forest}}\\ 
        \hline
          \textbf{Accuracy}& \textbf{Train} & \textbf{Test} & \textbf{Train} & \textbf{Test}\\
        \hline
          Overall & 84.2\% & 84.1\% & 89.1\% &  86.4\% \\
          LargeLabel & 91.9\% & 90.7\% & 99.7\% & 93.6\% \\
        \specialrule{1pt}{0pt}{0pt}
          \textbf{Time Factor} &  &  &  &  \\
          \hline
          Predicted & 0.914
          & 0.919
          & 0.884
          & 0.910
          \\
          Virtual Best & 0.862 & 0.863 & 0.862 & 0.863 \\
        \hline
\end{tabular}
\caption{Comparison of accuracy and runtime factors on training and test sets for final linear and random forest models}
\label{tab:MainResult}
\end{minipage}
\hspace{0.1\textwidth}
\begin{minipage}{0.4\textwidth}
\centering
\begin{tabular}{|c|c|c|}
    \hline
    & \textbf{Linear} & \textbf{Forest}\\ 
        \hline
        \textbf{Accuracy} &  & \\
        \hline
           Overall & 80.8\% & 85.2\%  \\
           LargeLabel & 77.0\% & 82.2\%  \\
        \specialrule{1pt}{0pt}{0pt}
          \textbf{Time Factor} &&\\
          \hline
          Predicted & 0.967 & 0.955  \\
          Virtual Best & 0.866 & 0.866  \\
        \hline
    \end{tabular}
    \caption{Results when using models trained on Xpress 9.6 for predicting 9.8 performance}
    \label{tab:9698}
\end{minipage}
\end{table}
\vspace*{-2.2\baselineskip}



\section{Conclusion}

In this paper, we investigated whether a global MINLP solver should always prioritize branching on fractional integer variables or whether allowing spatial branching earlier can lead to faster overall performance. Using a heterogeneous dataset of public and industrial MINLP instances, we trained linear and random forest regression models to predict the relative performance of two established branching rules. 
Our experiments demonstrated that linear models can achieve an 8\% reduction in mean solving time. We further showed that the learned models remain effective across solver versions, indicating robustness to underlying algorithmic changes.

A natural next step is to implement the regressor directly inside the FICO\textsuperscript{®} Xpress Global solver to validate its impact in a production setting and to extend them to other solvers, like the SCIP open-source MINLP solver. Finally, while our models select a single branching rule for the entire branch-and-bound tree, a more fine-grained approach, such as dynamically choosing rules for different phases of the solve~\cite{berthold2018feasibility} or adaptively at individual nodes of the branch-and-bound tree, represents a promising direction for further research.

\begin{credits}
\subsubsection{\discintname}
Timo Berthold is an employee of FICO.
\subsubsection{Acknowledgements.}
We thank Tristan Gally for his support with Xpress implementations and Ksenia Bestuzheva and Stefan Vigerske for the valuable discussions on SCIP features.
The work for this article was supported through the Research Campus
Modal funded by the German Federal Ministry of Education and Research (fund numbers 05M14ZAM,05M20ZBM).
\end{credits}

\bibliographystyle{splncs04} 
\bibliography{masterrefs} 

\end{document}